\documentclass[12pt,a4paper]{amsart} 
\usepackage{graphicx,latexsym,amssymb,color}
\theoremstyle{plain} 
\newtheorem{thm}{Theorem}[section]

\newtheorem{cor}[thm]{Corollary} 
\theoremstyle{definition}
 
\theoremstyle{remark} 
\newtheorem{rem}[thm]{\rm{Remark}} 
 
\numberwithin{equation}{section}
\numberwithin{figure}{section}

\begin{document} 
\title[Whitehead double and Milnor invariants]{Whitehead double and Milnor invariants} 
\author[J.B. Meilhan]{Jean-Baptiste Meilhan} 
\address{Institut Fourier\\
         Universit\'e Grenoble 1 \\
         100 rue des Maths - BP 74\\
         38402 St Martin d'H\`eres , France}
	 \email{jean-baptiste.meilhan@ujf-grenoble.fr}
\author[A. Yasuhara]{Akira Yasuhara} 
\address{Tokyo Gakugei University\\
         Department of Mathematics\\
         Koganeishi \\
         Tokyo 184-8501, Japan}
	 \email{yasuhara@u-gakugei.ac.jp}

\thanks{ 
The second author is partially supported by a Grant-in-Aid for Scientific Research (C) 
($\#$20540065) of the Japan Society for the Promotion of Science.}

%\date{\today}
%
\begin{abstract} 
We consider the operation of Whitehead double on a component of a link and study the behavior of 
Milnor invariants under this operation.  
We show that this operation turns a link whose Milnor invariants of length $\leq k$ are all 
zero into a link with vanishing Milnor invariants 
of length $\leq 2k+1$, and we provide formulae for the first non-vanishing ones.  
As a consequence, we obtain statements relating the notions of link-homotopy and self 
$\Delta$-equivalence via the Whitehead double operation.  
By using our result, we show that a Brunnian link $L$ is link-homotopic to the unlink if and only if 
a link $L$ with a single component Whitehead doubled is self $\Delta$-equivalent to the unlink. 
\end{abstract} 
\maketitle 
\section{Introduction}
In this paper, we consider the operation of Whitehead double, more generally of Whitehead $n$-double, 
on a component of a link, and we study the behavior of Milnor invariants under this operation.  
Milnor invariants $\overline{\mu}_L(I)$ of an $m$-component link $L$, where $I=i_1i_2...i_k$ with $1\le i_j\le m$, 
can be thought of as some sort of \lq\lq higher order linking number'' of the link.  
See Section~\ref{milnor} for a definition.  

A typical example is the Whitehead link, which is a Whitehead double of the Hopf link.  
The linking number of the Hopf link (which coincides with Milnor invariant $\overline{\mu}(12)$) is $\pm 1$, 
whereas the Whitehead link has linking number $0$.  
On the other hand, the Whitehead link has some nontrivial higher 
order Milnor invariants: its Sato-Levine invariant for instance, 
which is equal to $-\overline{\mu}(1122)$, is $\pm 1$.  
Our main result, stated below, generalizes this observation.  

\begin{figure}[!h]
\includegraphics{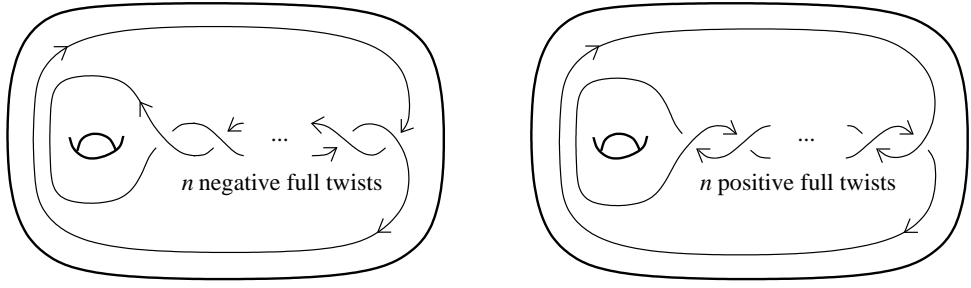}
\caption{The knot $\mathcal{W}_n$ for $n<0$ and $n>0$ respectively.} \label{wh}
\end{figure}

Let $K$ be a component of a link $L$ in $S^3$, regarded as $h(\{{\bf0}\}\times S^1)$ 
for some embedding $h:D^2\times S^1\rightarrow S^3\setminus(L\setminus K)$,   
such that $K$ and $h((0,1)\times S^1)$ have linking number zero.  
Let $n$ be a (nonzero) integer.  Consider in the solid torus $T=D^2\times S^1$ the knot 
$\mathcal{W}_n$ depicted in Figure~\ref{wh}.  
The knot $h(\mathcal{W}_n)$ is called the \emph{Whitehead $n$-double of $K$}, and it is denoted by $W_n(K)$.    

Given an $m$-component link $L=K_1\cup...\cup K_m$ in $S^3$, we denote by $W^i_n(L)$ the link 
$\left(L\setminus K_{i}\right)\cup W_n(K_{i})$ obtained by Whitehead $n$-double on the $i^{th}$ component of $L$.  
Note that the case $n=\pm1$ coincides with the usual notion of (positive or negative) Whitehead double.  

\begin{thm}\label{double}
Let $L$ be an $m$-component link in $S^3$, and let $n(\neq 0)$ be an integer.  
If all Milnor invariants $\overline{\mu}_L(Ji)$ of $L$ of length $|Ji|\leq k$ are zero $(k\ge 1)$, 
then all Milnor invariants $\overline{\mu}_{W^i_{n}(L)}(Ii)$ of $W^i_{n}(L)$ of
length $|Ii|\leq2k+1$ are zero.  
Moreover, if $\overline{\mu}_L(Pi)\ne 0$, $\overline{\mu}_L(Qi)\ne 0$ with $P=p_1 p_2 ... p_k,~Q=q_1 q_2 ... q_k$ 
$($possibly $P=Q)$ such that $p_j\ne i,~q_j\ne i$ for all $1\le j\le k$, 
then we have the following formulae for the first non-vanishing Milnor invariants of $W^i_{n}(L)$
           \begin{displaymath}
           \left\{ \begin{array}{l}
           \overline{\mu}_{W^i_{n}(L)}(PiQi)=2n\overline{\mu}_L(Pi)\overline{\mu}_L(Qi), \\
           \overline{\mu}_{W^i_{n}(L)}(PQii)=-n\overline{\mu}_L(Pi)\overline{\mu}_L(Qi).
           \end{array} \right.
           \end{displaymath}
\end{thm}

\begin{rem} \label{satolevine}
In the case of a $2$-component link, the formulae given in Theorem~\ref{double} for the 
first nonvanishing Milnor invariants of $W^i_n(L)$ provide, as an immediate corollary, 
a generalization of a result of Shibuya and the second author \cite{SY} as follows:  \\
Let $L=K_1\cup K_2$ in $S^3$.    
Let $n\ne 0$ be an integer, and let $W_n(L)$ be obtained by Whitehead $n$-double on a component of $L$.  
Then the Sato-Levine invariant $\beta_2$ of $W_n(L)$ satisfies 
 \[ \beta_2(W_n(L)) = n\left(lk(K_1,K_2)\right)^2. \]
\noindent (Note that the Sato-Levine invariant of $W_n(L)$ is well-defined, 
as Theorem~\ref{double} ensures that the link has zero linking number.)
\end{rem}

\begin{figure}[!h]
\includegraphics{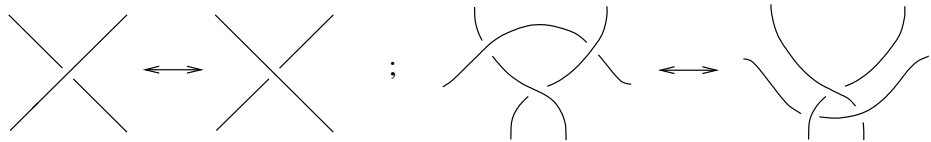}
\caption{A crossing change and a $\Delta$-move} \label{delta}
\end{figure}

Recall that two links are \emph{link-homotopic} if they are related by a 
sequence of ambient isotopies and \emph{self crossing changes}, 
which are crossing changes involving two 
strands of the same component, see the left-hand side of Figure~\ref{delta}.  
In particular, a link is called \emph{link-homotopically trivial} if it is link-homotopic to the unlink.  
It has long been known that Milnor invariants with no repeating indices 
are invariants of link-homotopy \cite{Milnor2}.  
Like crossing change, the $\Delta$-move is an unknotting operation \cite{MN}.  
Here we consider the notion of \emph{self $\Delta$-move} for links, 
which is a local move as illustrated in the right-hand side of Figure~\ref{delta} 
involving three strands of the same component. 
Two links are {\em self $\Delta$-equivalent} if they are related  
by a finite sequence of ambient isotopies and self $\Delta$-moves. 
Self $\Delta$-equivalence is a generalized link-homotopy, i.e., 
self $\Delta$-equivalence implies link-homotopy. 
The self $\Delta$-equivalence was introduced by Shibuya \cite{Shi,Shi1}, and was 
subsequently studied by various authors 
\cite{FY2,NO,NS,NSY,SY2,SY,yasuhara}.
A link is \emph{self $\Delta$-trivial} if it is self $\Delta$-equivalent to 
the unlink.

The following is a consequence of our main result.  

\begin{cor} \label{self}
Let $L$ be an $m$-component link in $S^3$ which is not link-homotopically trivial.  
Then, for any $n(\neq 0)$ and  $i~(1\le i\le m)$,   
$W^i_n(L)$ is not self $\Delta$-trivial.  
\end{cor}

Recall now that a link $L$ is Brunnian if all (proper) sublinks of $L$ are trivial.  
The next result shows that the converse of Corollary~\ref{self} also holds for Brunnian links.  
\begin{thm} \label{br}
Let $L$ be an $m$-component Brunnian link in $S^3$. Let $n(\neq 0)$ and $i~(1\le i\le m)$ be integers.  
Then $L$ is link-homotopically trivial if and only if $W^i_n(L)$ is self $\Delta$-trivial.  
\end{thm}

Observe that an $m$-component Brunnian link always has vanishing Milnor invariants of length $\le m-1$ 
since these are Milnor invariants of sublinks of a Brunnian link, which are trivial links.   
So Theorem~\ref{double} implies that all Milnor invariants of $W^i_{n}(L)$ of length $\le 2m-1$ are zero 
for any choice of $1\le i\le m$ and $n(\neq 0)$.  
In other words, for $m$-component Brunnian links, Whitehead doubling kills all Milnor invariants of length $\le 2m-1$.  
It follows from a more general result (stated and proved in Section~\ref{sat}) that an additional Whitehead doubling, 
on either the same or another component of the link, actually kills \emph{all} Milnor invariants, as the resulting 
link is always a boundary link, see Corollary~\ref{boundary}.     

The rest of the paper is organized as follows.  
In Section~\ref{milnor} we recall the definition of Milnor invariants and 
prove Theorem~\ref{double}.  % and Corollary~\ref{self}.  
In Section~\ref{sde} we prove the two statements relating Whitehead doubling 
and self $\Delta$-equivalence, 
namely Corollary~\ref{self} and Theorem~\ref{br}.  
In Section~\ref{sat} we consider more general satellite constructions, 
involving a knot which is null-homologous in the solid torus.  
When applied twice to a Brunnian link, such a construction always yields 
a boundary link.  
\section{Milnor invariants} \label{milnor}
%
%
%
%\subsection{Definition of Milnor Invariants} \label{mil}

J. Milnor defined in \cite{Milnor,Milnor2} a family of invariants of oriented, ordered 
links in $S^3$, known as Milnor's $\overline{\mu}$-invariants. 

Given an $m$-component link $L$ in $S^3$, denote by $\pi(L)$ the 
fundamental group of $S^3\setminus L$, and by $\pi_q(L)$ 
the $q^{th}$ subgroup of the lower central series of $\pi(L)$.  
We have a presentation  of $\pi(L)/ \pi_q(L)$ with $m$ generators, 
given by a meridian $\alpha_i$ of the $i^{th}$ component of $L$.  
So for $1\le i\le m$, the longitude $l_i$ of the $i^{th}$ component of $L$ 
is expressed modulo $\pi_q(L)$ as a word  
in the $\alpha_i$'s (abusing notations, we still denote this word by $l_i$).  

The \emph{Magnus expansion} $E(l_i)$ of $l_i$ is the formal power series in 
non-commuting variables $X_1,...,X_m$ obtained by 
substituting $1+X_j$ for $\alpha_j$ and $1-X_j+X_j^2-X_j^3+\cdots$ for $\alpha_j^{-1}$, $1\le j\le m$.  

Let $I=i_1 i_2 ...i_{k-1} j$ be a multi-index (i.e., a sequence of possibly repeating 
indices) among $\{1,...,m\}$. 
Denote by $\mu_L(I)$ the coefficient of $X_{i_1}\cdots X_{i_{k-1}}$ in the Magnus expansion $E(l_j)$.  
\emph{Milnor invariant} $\overline{\mu}_L(I)$ is the residue class of 
$\mu_L(I)$ modulo the greatest common divisor of 
all $\mu_L(J)$ such that $J$ is obtained from $I$ by removing at least one index, 
and permutating the remaining indices cyclically.  
We call $|I|=k$ the \emph{length} of Milnor invariant $\overline{\mu}_L(I)$.  
 
The indeterminacy comes from the choice of the meridians $\alpha_i$ or, 
equivalently, from the indeterminacy of 
representing the link as the closure of a string link \cite{HL}.  
\begin{proof}[{\rm Proof of Theorem~\ref{double}}] %\label{W}
%In this section we prove Theorem~\ref{double}.  
Without loss of generality, we may suppose that $i=m$. 
We give the proof of the case $n<0$.  The case $n>0$ is strictly similar and 
we omit it.  

We denote by $\alpha_1$,...,$\alpha_{m-1}$, $\alpha_m$ and $a$ 
meridians of $K_1$,...,$K_{m-1}$, $K_m$ and $W_{n}(K_m)$ respectively, 
such that 
$\alpha_1,...,\alpha_{m}$ generate $\pi(L)/\pi_q(L)$ and 
$\alpha_1,...,\alpha_{m-1},a$ generate 
$\pi(W^m_{n}(L))/\pi_q(W^m_{n}(L))$.  

The Magnus expansion of the longitude $l_m\in \pi(L)/\pi_q(L)$ of $K_m$, 
written as a word in $\alpha_1,...,\alpha_m$, has the form 
\[ E(l_m)=1+\sum{\mu}_L(i_1...i_r,m)X_{i_1}...X_{i_r} = 1+f(X_1,...,X_m), \]
where $E(\alpha_i)=1+X_i$ for all $1\le i\le m$.   

Now consider the Whitehead $n$-double of $K_m$, and consider $2n+1$ elements 
$a_0$, $a_1$, ... , $a_{2n}$ of $S^3\setminus W^m_{n}(L)$ 
as represented in Figure~\ref{W1}. Let $\phi(l_m)=l$, where $\phi: \pi(L)/\pi_q(L) \rightarrow \pi(W^m_{n}(L))/\pi_q(W^m_{n}(L))$ is the natural map that maps 
$\alpha_i$ to itself ($1\le i\le m-1$) and maps $\alpha_m$ to
$a_{2n}^{-1}a$. (Abusing notation, we still denote by $a_i$, $0\le i\le 2n$, 
the corresponding elements in $\pi(W^m_{n}(L))/\pi_q(W^m_{n}(L))$.)

\begin{figure}[!h]
\input{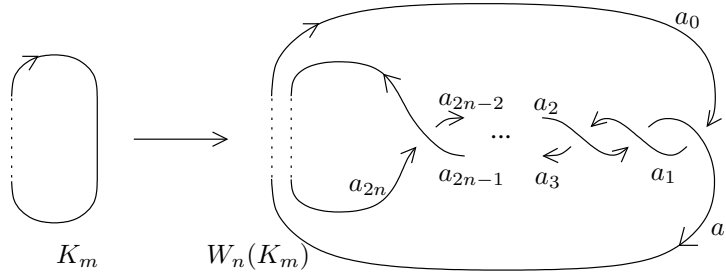} 
\caption{The Whitehead $n$-double of $K_m$ for $n<0$.} \label{W1}
\end{figure}

It follows from repeated uses of Wirtinger relations that 
\begin{displaymath}
\left\{ \begin{array}{ll}
a_0=l^{-1}al, & \\
a_{2r}=R^raR^{-r}, & \textrm{for all $r\ge 1$} \\
a_{2r+1}=R^raR^{-(r+1)}, & \textrm{for all $r\ge 0$} 
\end{array} \right.
\end{displaymath}
\noindent where $R=al^{-1}a^{-1}l$.  
In particular we have that  
   $$ \phi(\alpha_m)=a_{2n}^{-1}a= R^n a^{-1}R^{-n}a.  $$
Let $E(a)=1+X$ denote the Magnus expansion of $a$.  
Observe that 
 \[\begin{array}{rrcl}
 & E(R)=E(al^{-1}a^{-1}l)&=&(1+X)E(l^{-1})(1-X)E(l)+{\mathcal O}_X(2)\\
 & &=&1+X-E(l^{-1})XE(l)+{\mathcal O}_X(2),
\\
\textrm{and } & & & \\
 & E(R^{-1})=E(l^{-1}ala^{-1})&=&E(l^{-1})(1+X)E(l)(1-X)+{\mathcal O}_X(2)\\
 & &=&1-X+E(l^{-1})XE(l)+{\mathcal O}_X(2),
\end{array}\]
where ${\mathcal O}_X(2)$ denotes terms which contain $X$ at least $2$ times.  So we have 

\begin{eqnarray*}
E(\phi(\alpha_m))&=&(1+X-E(l^{-1})XE(l))^n(1-X)\\
&&\times(1-X+E(l^{-1})XE(l))^n(1+X)+{\mathcal O}_X(2)\\
&=&(1+nX-nE(l^{-1})XE(l))(1-X)\\
&&\times(1-nX+nE(l^{-1})XE(l))(1+X)+{\mathcal O}_X(2),\\
&=&1+{\mathcal O}_X(2).
\end{eqnarray*}
This implies that  
   $$ E(l)=1+f(X_1,...,X_{m-1},{\mathcal O}_X(2))=1+f_1(X_1,...,X_{m-1})+f_2(X_1,...,X_{m-1},X), $$
\noindent where 
$$f_1(X_1,...,X_{m-1})=f(X_1,...,X_{m-1},0)\in {\mathcal O}(k)$$
 and 
$$f_2(X_1,...,X_{m-1},X)=f(X_1,...,X_{m-1},{\mathcal O}_X(2))-f_1(X_1,...,X_{m-1})
\in {\mathcal O}(k+1),$$ 
and ${\mathcal O}(u)$ denotes terms of degree at least $u$ 
(the degree of a monomial in the $X_j$ is simply defined by the sum of the powers).  
Similarly we have 
   $$ E(l^{-1})=1+g(X_1,...,X_{m-1},{\mathcal O}_X(2))=1+g_1(X_1,...,X_{m-1})+g_2(X_1,...,X_{m-1},X), $$
\noindent where $g_1(X_1,...,X_{m-1})\in {\mathcal O}(k)$ and $g_2(X_1,...,X_{m-1},X)\in {\mathcal O}(k+1)$.  

Let $f_1,f_2,~g_1,~g_2$ denote $f_1(X_1,...,X_{m-1})$, $f_2(X_1,...,X_{m-1},X)$, 
$g_1(X_1,...,X_{m-1})$, $g_2(X_1,...,X_{m-1},X)$ respectively, 
and set $f=f_1+f_2$ and $g=g_1+g_2$. 
Set $E(a^{-1})=1-X+X^2-X^3+\cdots=1+Y$. 
Note that $(1+f)(1+g)=(1+g)(1+f)=1$ and $(1+X)(1+Y)=(1+Y)(1+X)=1$, 
hence $f+g=-fg=-gf\in {\mathcal O}(2k)$ and $X+Y=-XY=-YX$.  
One can check, by induction, that 
$$
\left\{ \begin{array}{l} 
E(R^n)=1+n(gY-Xf+XgY+gYf)+{\mathcal O}(2k+2), \\
E(R^{-n})=1+n(Xf-gY+XfY+gXf)+{\mathcal O}(2k+2), \\
E((a^{-1}R)^n)=(1+Y)^n+(1+Y)^nf - f(1+Y)^n+n(gYf- fgY)+{\mathcal O}(2k+2).
\end{array} \right.
$$

Since the preferred longitude $L_m$ of $W^m_{n}(K_m)$ is presented in $\pi(W^m_{n}(L))/\pi_q(W^m_{n}(L))$ 
by the word 
 $$ L_m = la^{-1}a_2^{-1}...a_{2n-2}^{-1}l^{-1}a_{2n-1}^{-1}a_{2n-3}^{-1}a_{3}^{-1}a_{1}^{-1}a^{2n} 
        = l(a^{-1}R)^{n}R^{-n}l^{-1}R^{n}a^{n},  $$
we have 
\begin{eqnarray*}
E(L_m)&=&(1+f)\big[(1+Y)^n+(1+Y)^nf - f(1+Y)^n+n(gYf- fgY)\big]\\
& &\times\big[1+n(Xf-gY+XfY+gXf) \big](1+g)\\%\big[1+n(gY-Xf+XgY+gYf) \big](1+X)^n\\
& &\times\big[1+n(gY-Xf+XgY+gYf) \big](1+X)^n\\%+{\mathcal O}(2k+2)\\
&=&\big[(1+Y)^n + n(2fXf - f^2X - Xf^2)\big](1+X)^n+{\mathcal O}(2k+2)\\
&=&1+n\left(2fXf - ffX - Xff\right)+{\mathcal O}(2k+2)
\end{eqnarray*}

Because $f\in {\mathcal O}(k)$, the first non-trivial terms in 
the Magnus expansion $E(L_m)$ are of degree $2k+1$.  
It follows that all Milnor invariants $\overline{\mu}_{W^m_{n}(L)}(Im)$ of length $|Im|\leq 2k+1$ of $W^m_{n}(L)$ are zero. 

Moreover, we actually have 
$$
E(L_m)=1+n\left(2f_1Xf_1 - f_1f_1X - Xf_1f_1\right)+{\mathcal O}(2k+2).
$$
So if $\overline{\mu}_L(Pm)\ne 0$, $\overline{\mu}_L(Qm)\ne 0$ for some multi-indices $P=p_1...p_k$, 
$Q=q_1...q_k~(P\ne Q)$ with $p_j\ne m$, $q_j\ne m$ for all $1\le j\le k$, 
then 
 $$f_1=\overline{\mu}_L(Pm)X_{p_1}...X_{p_k}+\overline{\mu}_L(Qm)X_{q_1}...X_{q_k}+{\mathcal O}(k), $$
and it follows from the above formula that 
\begin{eqnarray*}
E(L_m)&=&1+2n\overline{\mu}_L(Pm)\overline{\mu}_L(Pm) X_{p_1}\cdots X_{p_k}XX_{p_1}\cdots X_{p_k}\\
&&+2n\overline{\mu}_L(Pm)\overline{\mu}_L(Qm) X_{p_1}\cdots X_{p_k}XX_{q_1}\cdots X_{q_k} \\
&&+2n\overline{\mu}_L(Qm)\overline{\mu}_L(Pm) X_{q_1}\cdots X_{q_k}XX_{p_1}\cdots X_{p_k}\\
&&+2n\overline{\mu}_L(Qm)\overline{\mu}_L(Qm) X_{q_1}\cdots X_{q_k}XX_{q_1}\cdots X_{q_k}\\
&&-n\overline{\mu}_L(Pm)\overline{\mu}_L(Pm) X_{p_1}\cdots X_{p_k}X_{p_1}\cdots X_{p_k}X \\
&&-n\overline{\mu}_L(Pm)\overline{\mu}_L(Qm) X_{p_1}\cdots X_{p_k}X_{q_1}\cdots X_{q_k}X \\
&&-n\overline{\mu}_L(Qm)\overline{\mu}_L(Pm) X_{q_1}\cdots X_{q_k}X_{p_1}\cdots X_{p_k}X \\
&&-n\overline{\mu}_L(Qm)\overline{\mu}_L(Qm) X_{q_1}\cdots X_{q_k}X_{q_1}\cdots X_{q_k}X \\
&&-n\overline{\mu}_L(Pm)\overline{\mu}_L(Pm) XX_{p_1}\cdots X_{p_k}X_{p_1}\cdots X_{p_k} \\
&&-n\overline{\mu}_L(Pm)\overline{\mu}_L(Qm) XX_{p_1}\cdots X_{p_k}X_{q_1}\cdots X_{q_k} \\
&&-n\overline{\mu}_L(Qm)\overline{\mu}_L(Pm) XX_{q_1}\cdots X_{q_k}X_{p_1}\cdots X_{p_k}\\
&&-n\overline{\mu}_L(Qm)\overline{\mu}_L(Qm) XX_{q_1}\cdots X_{q_k}X_{q_1}\cdots X_{q_k} 
+ {\mathcal O}(2k+1)
\end{eqnarray*}
which implies the desired formulae for the first nonvanishing Milnor invariants of $W^m_{n}(L)$.  
\end{proof}

\begin{rem}
One may wonder what happens when we consider, in the definition of a 
Whitehead $n$-double, an odd number $2p+1$ of half-twists 
in place of $n$ full twists.  For a link $L$, denote by $W^i_{odd}(L)$ any 
link obtained by such a satellite construction with an odd number of 
half-twists on the $i^{th}$ component of $L$.  Then 
we can prove the following: 
If all Milnor invariants of $L$ with length $\le k$ vanish, then 
for any multi-index $Ii$ with $|Ii|\le k+1$, 
$\overline{\mu}_{W^i_{odd}(L)}(Ii)=2^{r_i+1}\overline{\mu}_L(Ii)$, where 
$r_i$ is the number of times that the index $i$ appears in $I$.    
\end{rem}
\section{On self $\Delta$-equivalence} \label{sde} 
In this section we provide the proofs for Corollary~\ref{self} and Theorem~\ref{br}.  

We need the following additional notation.  
Given a multi-index $I$, we denote by $r(I)$ the maximum number of times that any index 
appears in $I$.  
For example, $r(1123)=2$ and $r(1233212)=3$.    
\begin{proof}[{\rm Proof of Corollary~\ref{self}}]
Let $L$ be an $m$-component link which is not link-homotopi\-cal\-ly trivial.  
Then by \cite{Milnor} there exists some multi-index $I=i_1...i_p$ with $r(I)=1$  
such that $\overline{\mu}_L(I)\ne 0$ and $\overline{\mu}_L(J)=0$ for all multi-index $J$ with length $|J|<|I|$ 
and $r(J)=1$. 

Let $n(\neq 0)$ and $i~(1\leq i\leq m)$ be integers. If $I$ does not contain $i$, then 
$\overline{\mu}_{W^i_n(L)}(I)=\overline{\mu}_L(I)\neq 0$. 
So $W^i_n(L)$ is not link-homotopically trivial. Hence  $W^i_n(L)$ is 
not self $\Delta$-trivial. Suppose that $I$ contains $i$. 
By \lq\lq cyclic symmetry'' (\cite[Theorem 6]{Milnor2}), we may assume that $i_p=i$. 
By Theorem~\ref{double}, the link 
$W^i_n(L)$ thus satisfies $\overline{\mu}_{W^i_n(L)}(M)\neq 0$ for some multi-index $M$ with $r(M)\le 2$.  
Since Milnor invariants with $r\leq 2$ are self $\Delta$-equivalence invariants \cite{FY}, 
$W^i_n(L)$ is not self $\Delta$-trivial. 
\end{proof}
\begin{proof}[{\rm Proof of Theorem~\ref{br}}] \label{prbr}
Let $L$ be an $m$-component Brunnian link.  Let $n(\neq 0)$ and $i~(1\le i\le m)$ 
be integers.  
By Corollary~\ref{self} we already know that $L$ is link-homotopically trivial 
if $W^i_n(L)$ is self $\Delta$-trivial.  
Let us prove that the converse is also true.  

The link $L$ being Brunnian, $\overline{\mu}_L(I)=0$ if $I$ does not contain 
an index in $ \{ 1,...,m \} $. 
Moreover, if $L$ is link-homotopically trivial, then $\overline{\mu}_L(I)=0$ 
for any $I$ with $r(I)=1$.  
In particular $\overline{\mu}_L(I)=0$ for all $|I|\le m$, and by Theorem~\ref{double} 
the link $W^i_n(L)$ thus satisfies $\overline{\mu}_{W^i_n(L)}(I)=0$ for all $|I|\le 2m+1$.  
This implies that 
$\overline{\mu}_{W^i_n(L)}(I)=0$ for any multi-index $I$ with $r(I)\le 2$.  
By \cite[Corollary 1.5]{yasuhara}, 
we have that $W^i_n(L)$ is self $\Delta$-trivial.  
\end{proof}
\section{From Brunnian links to boundary links} \label{sat}
\subsection{Boundary links from satellite construction}
In this section we consider a more general satellite construction.  % than the Whitehead $n$-double.  

Let $L=K_1\cup...\cup K_m$ be an $m$-component link in $S^3$, and let 
$h_i:D^2\times S^1\rightarrow S^3$ be an embedding such that  
$h_i(\{{\bf 0}\}\times S^1)$ is the $i^{th}$ component $K_i$ of $L$ 
(as in the introduction, we assume that $K_i$ and 
$h((0,1)\times S^1)$ have linking number zero).  
Now, instead of the knot $\mathcal{W}_n$ depicted in Figure~\ref{wh}, 
consider in the solid torus $T=D^2\times S^1$ 
a fixed knot $\mathcal{K}$ which is null-homologous in $T$.  
Denote by $W^i_{\mathcal{K}}(L)$ the link $\left(L\setminus K_{i}\right)\cup h_i(\mathcal{K})$.  
We have the following result.
\begin{thm} \label{thsat}
Let $L=K_1\cup...\cup K_m$ be an $m$-component link in $S^3$, 
and let $\mathcal{K}$, $\mathcal{K}'$ be two null-homologous knots in the solid torus $T$.  
Then 
\begin{enumerate}
\item[(i)]  If $L\setminus K_i$ is a boundary link, then $W^i_{\mathcal{K}}(W^i_{\mathcal{K}'}(L))$ 
is a boundary link.
\item[(ii)] If $L\setminus (K_i\cup K_j)$ is a boundary link and $K_i\cup K_j$ is null-homotopic 
in $S^3\setminus\left(L\setminus (K_i\cup K_j)\right)$, 
            then $W^i_{\mathcal{K}}(W^j_{\mathcal{K}'}(L))$ is a boundary link. 
\end{enumerate}
\end{thm}

Note that in particular a Brunnian link $L$ always satisfies the conditions in (i) and (ii).  
It follows that a link obtained from a Brunnian link by taking 
twice Whitehead double (on either the same or another component of the link) 
kills \emph{all} Milnor invariants.   
\begin{cor} \label{boundary}
Let $L$ be an $m$-component Brunnian link in $S^3$. 
Let $p,q ~(pq\neq 0)$ and $i,j\in  \{ 1,...,m \} $ $($possibly equal$)$ be integers.  
Then the link $W^{i,j}_{p,q}(L)$, obtained by respectively Whitehead $p$-double 
and Whitehead $q$-double on the $i^{th}$ and 
$j^{th}$ components of $L$, is a boundary link.  
\end{cor} 
\noindent Figure~\ref{ex} below illustrates this result in the case of the Borromean rings $B$.  

\begin{figure}[!h]
\includegraphics{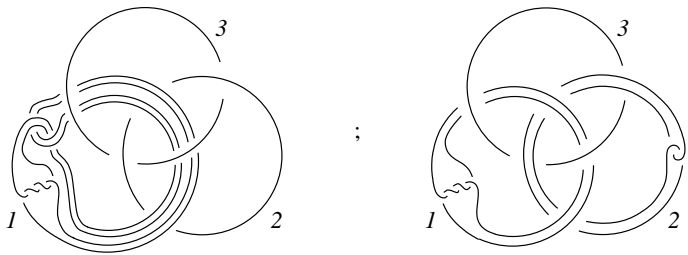}
\caption{The boundary links $W^{1,1}_{-4,2}(B)$ and $W^{1,2}_{-4,2}(B)$.} \label{ex}
\end{figure}

\subsection{Proof of Theorem~\ref{thsat}}
Before proving Theorem~\ref{thsat}, 
we will introduce the notion of band presentation of a link.  

\begin{figure}[!h]
\input{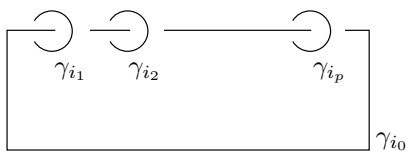} 
\caption{The link $L_i$.} \label{canonical}
\end{figure}

Let $L_i=\gamma_{i0}\cup \gamma_{i1}\cup \gamma_{i2}\cup... \cup \gamma_{ip_i}$ be a link as 
illustrated in Figure~\ref{canonical}.  
Let $L_1\cup...\cup L_m$ be a split union of the links $L_1$,...,$L_m$, and let 
$\Delta=\bigcup \Delta_{ij}$ be a disjoint union of disks $\Delta_{ij}$ 
($1\le i\le m$ ; $1\le j\le p_i$) such that $\partial \Delta_{ij}=\gamma_{ij}$ and 
$\Delta_{ij}\cap \left(\bigcup_{k} \gamma_{k0}\right)=\Delta_{ij}\cap \gamma_{i0}$ 
consists of a single point.  
It is known \cite{suzuki} that an $m$-component link $L$ in a $3$-manifold $M$ which is 
null-homotopic in $M$ can be expressed as a 
band sum of $L_1\cup...\cup L_m$, which is contained in a $3$-ball in $M$, 
along mutually disjoint bands $b_{ij}$ ($1\le i\le m$ ; $1\le j\le p_i$), 
disjoint from $\mathrm{int}\Delta$, 
such that $b_{ij}$ connect $\gamma_{ij}$ and 
$\left(\bigcup_{k} \gamma_{k0}\right)$. 
\footnote{The result is given in \cite{suzuki} for \emph{knots} in $S^3$, but 
it can be easily extended to the link case.} 
This presentation is called a \emph{band presentation} of $L$, and $L_1\cup...\cup L_m$ 
is called the \emph{base link}.  

\begin{proof}[{\rm Proof of Theorem~\ref{thsat}}]
{\rm (i)}~~~ We may suppose that $i=m$ without loss of generality.  
Since $K_1\cup...\cup K_{m-1}$ is a boundary link, it bounds a disjoint 
union of surfaces $E=E_1\cup...\cup E_{m-1}$.  
Denote by $W_{\mathcal{K}'}(K_m)$ the $m^{th}$ component of $W^m_{\mathcal{K}'}(L)$.  
Since $W_{\mathcal{K}'}(K_m)$ is null-homologous in $h_m(D^2\times S^1)$, it is 
null-homotopic in $S^3\setminus (L\setminus K_m)$. Hence there is a band presentation 
of $W_{\mathcal{K}'}(K_m)$ such that the base link is disjoint from $E$ and such that the 
intersections of each band and $E$ are ribbon singularities.  
So $W_{\mathcal{K}'}(K_m)\cap E$ is a union of copies of $S^0$, which are the 
endpoints of these ribbon singularities.  
By tubing the surfaces $E_i$ suitably at these endpoints, we obtain a union of mutually disjoint surfaces 
$F_1$,...,$F_{m-1}$ such that $F_i=\partial K_i$ and 
$F_i\cap W_{\mathcal{K}'}(K_m)=\emptyset$ for all $1\le i\le m-1$.  
Since the $m^{th}$ component of $W^m_{\mathcal{K}}(W^m_{\mathcal{K}'}(L))$ 
bounds a Seifert surface $F_m$ in a regular neighborhood of $W_{\mathcal{K}'}(K_m)$, it follows that 
the components of $W^m_{\mathcal{K}}(W^m_{\mathcal{K}'}(L))$ bound $m$ mutually disjoint 
Seifert surfaces $F_1$,...,$F_m$.  

{\rm (ii)}~~~ We may suppose that $i=m-1$ and $j=m$ without loss of generality.  
$K_1\cup...\cup K_{m-2}$ being a boundary link, it bounds a 
disjoint union of surfaces $E=E_1\cup...\cup E_{m-2}$.  
Since $K_{m-1}\cup K_m$ is null-homotopic in $S^3\setminus (K_1\cup...\cup K_{m-2})$, there is a 
band presentation of $K_{m-1}\cup K_m$  
such that the base link is disjoint from $E$ and such that the intersections of each band 
and $E$ are ribbon singularities.  
By tubing the surfaces $E_i$ suitably at the endpoints of theses singularities, 
we obtain a union of mutually disjoint surfaces 
$F_1$,...,$F_{m-2}$ such that $F_i=\partial K_i$ and $F_i\cap (K_{m-1}\cup K_m)=\emptyset$ 
for all $1\le i\le m-2$.  
Since the $(m-1)^{th}$ and $m^{th}$ components of $W^{m-1}_{\mathcal{K}}(W^m_{\mathcal{K}'}(L))$ 
bound a disjoint union $F_{m-1}\cup F_m$ of Seifert surfaces in a regular neighborhood of 
$K_{m-1}\cup K_m$, it follows that the components of $W^{m-1}_{\mathcal{K}}(W^m_{\mathcal{K}'}(L))$ 
bound $m$ mutually disjoint Seifert surfaces $F_1$,...,$F_m$.
\end{proof}

\end{document}